\documentclass[10pt]{article}

\usepackage[]{amsmath,amscd,amsfonts,amsthm,amssymb}
\usepackage[greek, english]{babel} 
\languageattribute{greek}{polutoniko}

\usepackage{graphicx}
\usepackage[all]{xy}
\input{diagxy}
\usepackage{psfrag}
\usepackage[lflt]{floatflt}
\newcommand{\bq}{\begin{quote}}
\newcommand{\eq}{\end{quote}}

\newcommand{\order}{\leqslant}

\newtheorem{Th}{Theorem}
\newtheorem{ax}{Axiom}
\newtheorem{lm}{Lemma}
\newtheorem{df}{Definition}
\newtheorem{pr}{Proposition}
\newtheorem{cl}{Corollary}
\newtheorem{re}{Remark}
\newtheorem{as}{Assumption}
\newtheorem{wg}{Wild Guess}
\newtheorem{ex}{Example}
\newcommand{\bth}{\begin{Th}\hspace{-5pt}{\bf .} \ }
\newcommand{\Eth}{\end{Th}}
\newcommand{\bax}{\begin{ax}\hspace{-5pt}{\bf .} \ }
\newcommand{\eax}{\end{ax}}
\newcommand{\blm}{\begin{lm}\hspace{-5pt}{\bf .} \ }
\newcommand{\elm}{\end{lm}}
\newcommand{\bdf}{\begin{df}\hspace{-5pt}{\bf .} \ }
\newcommand{\edf}{\end{df}}
\newcommand{\bpr}{\begin{pr}\hspace{-5pt}{\bf .} \ }
\newcommand{\epr}{\end{pr}}
\newcommand{\bcl}{\begin{cl}\hspace{-5pt}{\bf .} \ }
\newcommand{\ecl}{\end{cl}}
\newcommand{\bre}{\begin{re}\hspace{-5pt}{\bf .} \ }
\newcommand{\ere}{\end{re}}
\newcommand{\bas}{\begin{as}\hspace{-5pt}{\bf .} \ }
\newcommand{\eas}{\end{as}}
\newcommand{\bwg}{\begin{wg}\hspace{-5pt}{\bf .} \ }
\newcommand{\ewg}{\end{wg}}
\newcommand{\bex}{\begin{ex}\hspace{-5pt}{\bf .} \ }
\newcommand{\eex}{\end{ex}}

\newcommand{\bit}{\begin{itemize}}
\newcommand{\eit}{\end{itemize}\par\noindent}
\newcommand{\ben}{\begin{enumerate}}
\newcommand{\een}{\end{enumerate}\par\noindent}
\newcommand{\beq}{\begin{equation}}
\newcommand{\eeq}{\end{equation}}
\newcommand{\beqa}{\begin{eqnarray*}}
\newcommand{\eeqa}{\end{eqnarray*}\par\noindent}
\newcommand{\beqn}{\begin{eqnarray}}
\newcommand{\eeqn}{\end{eqnarray}\par\noindent}
\newcommand{\grk}{\selectlanguage{polutonikogreek}}
\newcommand{\eng}{\selectlanguage{english}} 
 
\def\nn{\mbox{I\hspace{-0.5mm}N}} % het getal kan aangepast worden naar smaak

\usepackage{fancybox}

\title{\huge On what Ontology Is and not-Is} 
\author{} 
\date{}

\begin{document}   

\maketitle  

\vspace{-2cm} 

%\centerline{\small \sc {[DRAFT]}}
\vspace{1cm} 

\centerline{Karin Verelst}
\vspace{1mm}
\par
{\scriptsize
\centerline{{\em FUND-CLEA}}
\vspace{1mm}
\centerline{{\em Vrije Universiteit Brussel}}
\centerline{{\em Pleinlaan 2, B-1050 Brussels}}
\centerline{{\em kverelst@vub.ac.be}}
} 

\bigskip
\bigskip
\bigskip

\begin{flushright}{\footnotesize {\em Indeed, language is an organ of perception, \\not simply a means of
communication.\\ {\em J. Jaynes},
Origin of Consciousness}}
\end{flushright}

\bigskip

\subsection*{\sc Introduction}

\frenchspacing In this contribution we place the relation between logic and classical metaphysics in a new light. We shall turn our attention to what Plato calls ``dialectics''  [{\em Sophist}, 253(c)], the method  (in the original sense of ``the way to follow towards'') by which true knowledge is acquired. This method remains fully embedded in his ontology, causing him, according to a widely held view, to lapse into elementary logical errors. A more subtle objection imputes on him the so-called ``Socratic fallacy''.\footnote{This attirbution of such fallacies started with P. Geach's influential article ``Plato's Euthyphro: An Analysis and Commentary'', {\em The Monist}, {\bf 50}, 1966, pp. 369-382.} We, on the contrary, hold that there are no such errors in Plato's method precisely {\em because} of its onto-logical nature: when studying Plato's ``logic'' one is simply studying his metaphysics from another angle; the structural properties of the latter can shed light on the logical structure of the former. Sayre puts it succinctly: {\em The connection between Forms and definitions is that participation in the Form constitutes what the formula defines}.\footnote{K.M. Sayre, {\em Plato's Late Ontology}, Princeton University Press, Princeton, New Jersey, 1983, p. 8.} This allows Plato to straightforwardly explain the nature of dialectics: {\em For it is from the mutual intertwinement of the Forms that reasoning comes forth} [{\em
Sophist}, 259(e)]. This is testified explicitly by Aristotle in [{\em Met. A}, 987b(31-32)]. It is not our intention to enter here into the delicate problem of the way things in the world participate into Forms in another ontological realm, but we accept  --- with Sayre --- that at least the later Plato was convinced to have found a satisfactorily solution to it. We can therefore both agree and disagree with M. Dixsaut in her penetrating study of the transformations Plato's notion of dialectics underwent throughout the dialogues:  {\em Que la dialectique ne soit pas une ``methode'', mais la science la plus haute ne faisant qu'un avec le chemin qu'elle se fraie est un trait essentiel de la philosophie Platonicienne, ins\'{e}parable \`{a} la fois de la repr\'{e}sentation qu'il se fait de la pens\'{e}e et de l'hypoth\`{e}se des formes.}\footnote{M. Dixsaut, {\em M\'{e}tamorphoses de la dialectique dans les dialogues de Platon}, Vrin, Paris, 2001, p. 9.} We comply to that, except for the fact that this is exactly why Plato's dialectics {\em is} a genuine method. When studying the classics, it is our conception of method that needs a change.\footnote{The idiosyncratic nature of Modern ideas on method is the subject of a paper entitled {\em Infinity and the Sublime}, which I presented at the Centre for Mathematical Sciences, Cambridge, October 2006 [forthcoming].} 

Our approach shall be to revisit Plato's idea that, in order to find out the specific Form of which a given thing is an instance, one has to proceed by dichotomic division. We construe this division as the logical backbone of Plato's late ontology, its ``pattern of reality'' according to Cherniss.\footnote{H. Cherniss, {\em The Riddle of the academy}, Russell \& Russell, New York, 1962, p. 42.} 
Following  a suggestion by G. Priest\footnote{G. Priest and R. Routley, ``First Historical Introduction. A Preliminary History of Paraconsistent and Dialethic Approaches'', in: {\em Paraconsistent Logic. Essays on the Inconsistent}, Philosophia Verlag, M\"{u}nchen, 1989, p. 19.} we shall then show that the logic shoring up Plato's metaphysical system ---  the doctrine of Ideas and the participation theory --- is {\bf paraconsistent}. As such it is the ``missing link'' between the paradoxical way of thinking of the pre-Socratics and classical onto-logic as established by Aristotle. Plato's teachings constitute the bridge that crosses the abyss between two utterly different forms of conscious relation to reality as it is ``given'': pre-Socratic world-articulation and classical (metaphysical)  world-description. In the first part of this paper I will sketch how precisely the ontologically paradoxical, {\em deictic} first person stand-point {\em in} reality of the archaic period differs from the essentially logical, third person view-point {\em on} reality as presented by classical metaphysics. In the second part I shall analyse how Plato tackles this problematic pre-Socratic legacy mainly in the form presented by Parmenides and his followers, and why he lays out the groundwork for almost all future attempts to deal with it: a {\em metaphysical system} in which the pre-Socratic paradoxes are resolved by {\bf a radical world division} in which {\em being} is separated from {\em non-being}, though in such a way that sufficient common ground remains present to fulfill the need for justification of the multiplicity and mutability we experience in `phaenomenal' reality. A metaphysical system is not just a grand story that explains the basic features of reality as a whole, nor even a characterisation of ``being {\em qua} being'', to quote the Peripatetic formula. It will be shown that Plato achieves this remarkable feat basically by transforming Zeno's paradoxical One-and-Many into a paraconsistent Large-and-Small. In the last part of this paper, a more in depth formal analysis is presented, and the link to Aristotle's system --- who ``completes'' the work both from the metaphysical and the logical point of view --- is shortly discussed.

\subsection* {\sc Pre-Socratic Paradoxicality}

 \noindent We shall at first briefly reconsider the often misunderstood
\vspace{-0.1 mm} pre-Socratic legacy, summed up concisely in the Parmenidean dictum\selectlanguage{polutonikogreek} t`o >e`on >'esti\selectlanguage{english} ({\em the Being-Now is} [{\small DK} {\scriptsize{28B 6}}]). In the present everything {\em
is}.\footnote{The fact that ``is'' with Parmenides cannot possibly refer to some abstract, eternal category --- as it does with Plato --- has been discussed by N.-L. Cordero, {\em Les deux chemins de Parm\'{e}nide}, Paris, 1984.} The reference to spatio-temporal {\em presence} is essential, although absent from the common translation {\em (the) Being is}. Restoring it turns Parmenides's claim from an  incomprehensible declaration with bizarre consequences into the infallibly true statement that it was intended to be.\footnote{K. Riezler, {\em Parmenides. Text,
\"{U}bersetzung, Einf\"{u}hrung und Interpretation}, Vittorio Klostermann, Frankfurt,
1970, p. 45-50. That the in origin dialectal difference between \grk t'o >'on \eng and \grk t`o >e`on \eng had
acquired philosophical significance becomes explicit in Diogenes of Appolonia, a
contemporary to Zeno, where \grk t`a >'onta \eng{\em [ta onta; the beings]}
are stable essences, while  \grk t`a >e`onta n~un  \eng {\em [ta
eonta nun; the beings-now]} are instable phaenomenological things. See L.
Couloubaritsis, {\em La Physique d'Aristote}, Ousia, Bruxelles, 1997, p. 308.} Indeed, you cannot {\em indicate} anything {\em here and now} that does {\em not} exist! I shall call the unmediated relation to reality constituted by this truly {\em indicative} capacity of a speaking ``world-centre'' which points at and thus points out while addressing another one with the same capacity {\em deictic}, in reference to and in honor of the French linguist Emile Benveniste, who introduced the notion of {\em deixis} in his investigations of the role of the ``first person'' in natural, spoken language.\footnote{E. Benveniste, ``Le langage et
l'exp\'{e}rience humaine'', in: {\em Probl\`{e}mes de linguistique  g\'{e}n\'{e}rale II},
Gallimard,  Paris, 1966, p. 69. I discussed this with respect to archaic reality-awarenessin detail in the third chapter of my PhD-thesis: K. Verelst, {\em De Ontologie van den Paradox, doctoral dissertation}, Vrije Universiteit Brussel, 2006.} Interestingly enough the situation in which a deicitc centre literally stands is ontologically paradoxical, for it instantiates a {\em coincidentia oppositorum}. At every moment one {\em is} and {\em is not} at the same place at the same moment, because the {\em now} can never be fixed or stabilised, so that identity does not exist outside this very same moment. This implies that the principle of identity does not hold. It moreover precedes spatio-temporal differentiation: I am Here-Now. It is not experienced as a problem {\em because} it is experienced, not merely thought. This is exactly what
Heraclitus pointed out with his famous riverfragment, which states that one cannot step into the same river twice  [{\small DK 22} {\scriptsize{B 49a}}].\footnote{One should appreciate the acuity of his disciple Cratylus who, seeing the consequences of the non-temporality of the now, states that one could not do it even {\em once}.} Revealingly enough, the later philosophical tradition saw him as the {\em opponent} of Parmenides.  I, however, take it that this deictical stance straightforwardly explains both Heraclitus and Parmenides, and thus that Conche is right when he says: {\em La pens\'{e}e de
Parm\'{e}nide se comprend d'une mani\`{e}re compl\`{e}te si l'on y voit une double
radicalisation de la pens\'{e}e d'H\'{e}raclite}.\footnote{M. Conche, {\em Parm\'{e}nide. Le Po\`{e}me: Fragments}, PUF, 1999, p. 26.} In this respect the pre-Socratics remain much closer to Homer than to
Aristotle, even when much closer to the latter in time: {\em Greeks like
Aristotle and we to-day have apparently attained to greater `detachment', power of
thinking in cold blood without bodily movement, as we have to a sharper discrimination
and definition of the aspects and phases of the mind's activity. It is with the
consciousness, the knowing self, the spectator aware of what happens within and without ($\ldots$) that a man would tend more particularly to identify himself}.\footnote{R.B. Onians, {\em The Origins of European Thought}, Cambridge University Press, Cambridge, 1994 [1951], p. 18.} The shift from first to third person induces the separation of subject and object, by loosening or even breaking the ties between a man's embodied awareness and his mental activity. In this sense one could say that, in the pre-Socratic era, the mental disposition prerequisite to the {\em possibility} of metaphysical thought was still lacking: the capacity to re-present the absent as present on a logically structured inner mental scenery, as described, incidentally, for the first time by Plato in his myth of the Cave, and internalised by Aristotle with his conception of {\em fantasia} as an inner mental screen [{\em De anima}, 3.3].\footnote{I owe this latter connection to F. De Haas,  Leiden University. For a discussion of Aristotle's notion of {\em fantasia}, see G. Watson, ``Fantasia in Aristotle, De Anima 3. 3'', {\em The Classical Quarterly}, New Series, Vol. 32, {\bf 1} 1982, pp. 100-113. A more recent study is A.D.R. Sheppard, ÒPhantasia and Inspiration in Neoplatonism.Ó In {\em Studies in Plato and the Platonic tradition}, M. Joyal (ed.), Aldershot, 1997, pp. 201-10.} Quoting von Fritz's comment on these developments {\em in extenso} is worthwile: \begin{quote} {\em \small (...) the concepts of the ``obscure" Heraclitus
are all perfectly clear and can be very exactly defined. In contrast, the empiricist Sextus,
whose arguments seem so clear and easy to many readers, has no clearly identifiable concept
of either {\em logos} [reason] or {\em nous} [mind] at all. {\em Nous} with Sextus is
either identified with {\em logos} or considered a manifestation of it. {\em Logos}, where
Sextus speaks in his own name, is most often ``logical reasoning" or the capacity of logical
reasoning (...). But where Sextus reports the views of other philosophers, logos becomes
just the alternative to {\em aisthesis} [perception], whatever this alternative may be,
and so loses all clearly identifiable meaning. {\bf \em Yet it is highly illustrative of the
change which the concept of nous had undergone} between Heraclitus and Sextus that Sextus,
in trying to explain Heraclitus' concept, begins by connecting it with a term the
preponderant meaning of which is ``reasoning" and ends by almost identifying it with
``sensual perception." {\bf \em Heraclitus' own concept of {\em nous}}, as we have seen, {\bf \em was clearly distinguished from both} but somewhat more nearly related to the latter than to
the former.} To the archaic mind, to think is to {\em see}, to have an insight, not to make an
inference from {\em a priori} praemises: {\em The fundamental meaning of the word \grk noe~in \eng in Homer is ``to realize or to understand a
situation''}.\footnote{K. Von Fritz, ``Nous, noein and their derivatives in pre-Socratic
philosophy''. Reprinted in: A. P. D. Mourelatos, {\em The pre-Socratics. A collection of
critical essays}, Princeton University Press, Princeton etc., 1992, pp. 42-43; p. 23 (my bold).}\end{quote} When one speaks, one does so not about a possible world, but about {\em (everything in) this world} --- \grk k'osmoc t'ode \eng according to
Heraclitus [{\small DK} {\scriptsize{22B 30}}] ---,  behind the surfaces of which nature likes to hide: {\em ($\ldots$) it is still the primary function of
the \grk n'ooc \eng to be in direct touch with ultimate reality}\footnote{K. Von Fritz, {\em o.c.}, p. 52.}, i.e., with {\em this world} as it is immediately and fully present and within deictic reach.\footnote{Conche translates ``ce
monde-ci''. M. Conche, {\em H\'{e}raclite. Fragments}, Presses Universitaires de France,
1986/1998, pp. 279-280.} The lack of temporality manifest in the deictic ``first person'' standpoint thus has a spatial counterpart. It is not spatiality that individuates a being, but a being-present that instantiates a space...\footnote{How this ``body-bound'' localisation of the of the spatiotemporal centre of reference amounts into a modern notion of relativity has been shown by H. Poincar\'{e}, ``La Relativit\'{e} de L'Espace'',  in  {\em Science et M\'{e}thode}, p. 92. I gave an overview of the issues involved in a presentation entitled {\em Deixis and Theoria.  Insight and Inference from Aristotle to Poincar\'{e}}, for the ECCO-seminar at the Vrije Universiteit Brussel, May 23,  2007. See also the remark in the DK apparatus, vol. I, p. 254.} ``Space'' nor ``time'' exist as independent backgrounds against wich a mentally representable event takes
place: {\em Il n'y a pas de cosmogonie [chez H\'{e}raclite] malgr\'{e}
l'apparence, parce qu'il n'y a pas de repr\'{e}sentation ($\ldots$)}.\footnote{J.
Bollack, H. Wismann, {\em H\'{e}raclite ou la
s\'{e}paration}, Ed. de Minuit, Paris, 1972, p. 49.} It is by indication that
the validity of an utterance is {\em shown}. But, as already said, this is a paradoxical situation, because this centered absolute
simultaneity excludes any rigid identity. 

The paradoxical ontology underlying these I-Here-Now utterances\footnote{I introduced these notions in my PhD-Thesis: K. Verelst, {\em De Ontologie van den Paradox}, doctoral dissertation, Vrije Universiteit Brussel, 2006.} is epitomised in the notorious paradoxes of Zeno, a disciple of Parmenides. If my deictic reading of the pre-Socratics is correct, then the key to a correct understanding of all --- both Plurality and Motion --- arguments based on Zeno's divisional procedure is that {\em they do not presuppose space, nor time.} Division merely requires extension (with a ``here'' and a ``there'') of a physical object present to the senses and takes
place simultaneously. Zeno shows (and this is the appropriate word) that plurality and motion are examples of the fundamental paradox of Being and non-Being as instantiated in the here-and-now. He nowhere {\em denies} the reality of plurality and change; he simply points out their paradoxical nature. This I call {\bf Zeno's deictic realism}.\footnote{Zeno's paradoxical procedure is analysed in detail from this perspective in: K. Verelst,  ``ZenoÕs Paradoxes. A Cardinal Problem. I. On Zenonian Plurality'', in: Paradox: Logical, Cognitive and Communicative Aspects.  Proceedings of the First International Symposium of Cognition, Logic and Communication, Series: The Baltic International Yearbook of Cognition, Logic and Communication, Vol. 1, University of Latvia Press, Riga, 2006.} One could say that his paradoxes are not problems but descriptions of real states of affairs. It is plain that Zeno's arguments are not a {\em reductio}, because the logical prejudice that something which implies paradoxes cannot ``really be there'' is still unthinkable, for logic itself does not yet exist!\footnote{Vastos's reading is interesting as a witness to our distorted understanding of this fact: {\em Three
distinct inferential sequences are joined to form a unified argument, exhibiting,
probably for the first time in a philosophical context, the {\em reductio} in its
most powerful form, ``if $P$, then $C$ and not-$C$. [Therefore not-$P$.]''} But the latter assertion is nowhere in the material. G. Vlastos, {\em Studies in Greek Philosophy, volume I: The
Presocratics}, Princeton University Press, Princeton, 1993, pp. 219-240.]}  Thus Simplicius's attestation that {\em [In his book, in which many
arguments are put forward,] he shows {\bf in each} that stating a plurality comes
down to stating a contradiction}: \grk kaj' <'ekaston de'iknusi, <'oti t~w poll`a e~>inai l'egonti sumb'anei tˆ >enant'ia l'egein \eng [Simpl., {\em Phys.}, 139 (5) ({\small DK}
{\scriptsize{29B 2}})] turns out to be entirely correct. 
In what are
traditionally called the plurality arguments, this contradiction appears in the
stature of the \grk mikr`a te e~>inai ka`i m'egala \eng --- the large[s]-and-small[s] [{\small DK}
{\scriptsize{29B 1}}] ---, the infinity of segments with and the infinity of segments without magnitude that result from Zeno's divisional procedure {\em simultaneously}.\footnote{The {\em te kai} construction \vspace{-0.5 mm} stresses the fact that both are the case at the
same time.} Their simultaneous presence is precisely what keeps an extended body {\em finite}. Although the argument is by no means a dilemma\footnote{Of the disjunctive particles \grk >'h? >'htoi \eng {\em [or]} needed to construe the argument as a dilemma there is no trace.}, one could use the dilemma-construction to prove this by contraposition.

Only after the invention of metaphysics - i.e., after Plato - it becomes
possible to understand the two key figures of pre-Socratic thought, Heraclitus and Parmenides, as standing in
contradiction to each other, as defending opposite ``worldviews". The supposedly contradicting conclusions deriving from pre-Socratic philosophy were of a major concern
to both Plato and Aristotle, because they challenged the existence of truth and certainty {\em about}
the world, by now conceived of as external, and therefore about the actions of human beings {\em in} it. One will appreciate the depth of the abyss that separes this abstracted, mediated third person relation to the world from the first person immediacy of the pre-Socratic reality-participation gone before. The uncertainty flowing forth from the attempt to {\em define} everything from the outside had given rise to a sceptical discipline, Sophism, that simply denied any relation between reality and what we can say about it.  Its subjectivism
stems from a radical empiricism, which holds that things are for me as I perceive them. The Sophists thus seem to speak in terms of the archaic {\em coincidentia}, but {no deictic link between perception and perceived is preserved, so that ``I'' can say anything about whatever: first person perceptive immediacy has been replaced by third person subjective empiricism. Nevertheless reminiscences to the old ontology lurk underneath, as is clear from Protagoras's famous quote {\em man is the measure of all things, {\bf \em existing and non-existing}}. Plato cites it [{\em Theaetetus} 160(d,e))] in an attempt to point out both its origins and consequences: since reality as we perceive it is always in a process of permanent change it implies also the non-existence of stable, individual things in the world. What follows for knowledge out of the ontology of permanent change is the ``everything is equally true", which reveals itself as {\bf the epistemological face of the coincidence of opposites}, according to Plato the very foundation of the Heraclitean universal change theory. Therefore Sophism is nothing else but an instance of the ontology of permanent change, already advanced --- according to Plato --- by earlier thinkers like Heraclitus. Now on this soil the principle of contradiction rests, because if you {\em allow} contradiction, you will be allowed to say whatever: {\bf \em ex falso quodlibet} ({\em Theaetetus}, 150, 182(e)-183(a)). The dramatic consequences of the shift in perspective from the first to the  third person could hardly be more clearly exposed: what was once the seat of deictic certainty had been transformed into the source of irreparable deception. But how to conceal the quest for certainty with our experience of permanent change
in {\em this world}? This will only be possible by stabilising human world-experience in a
world-{\em picture}, strong enough to survive the paradoxical present into the past and the
future.\footnote{The necessary condition that made possible this construction of stabilising
world-pictures or ``worldviews", was the earlier coming-to-be of the ``inner mind-space", in
which the non-present could be re-presented as present. This crucial notion is introduced and discussed in J. Jaynes, {\em The origin of consciousness in the breakdown of the bicameral mind}, Houghton Mifflin Company, Boston, 1976, , p. 54 sq.).} Strong
foundations must be laid to grant the possibility to experience entities as objects outside
of the stream of events, and therefore to speak about them in an ``objective" way. Thence the well-known scaffold of Plato's solution: the division
of the world in two separated, though connected, layers with a different ontological status: a basic unchangeable, motionless `Parmenidean' or `Eleatic' {\em Being} which grants certainty about both objects and names, and an unstable `Heraclitean' {\em non-Being}, which allows for the change and motion in the world as presented to our senses. The paradox of the being of non-being is resolved by Plato and Aristotle by taking non-Being as a modus of Being -- of a lesser ontological degree.\footnote{Pelletier gives a detailed overview of the different possible positions, with reference to certain passages in the {\em Sophist} discussing what he calls ``Parmenids' problem'; F.J. Pelletier, {\em Parmenides, Plato, and the Semantics of Not-Being}, Univ. of chicago Press, Chicago and London, chapters two \& four.} This feature of a two-layered world (a ``world behind the
world"), is precisely what makes their worldviews {\bf metaphysical}. The widely acclaimed chasm between Plato's and Aristotle's systems can be understood from this point of view as a mere variation on a common ontological theme required to solve the paradoxes. One can charcterise the metaphysical variants by the paradoxes they presuppose: while Plato in his dialogues discusses mainly Zeno's paradoxes of plurality, Aristotle confines himself almost
exclusively to the paradoxes of motion. This amounts in their `static' vs. `dynamic' metaphysics respectively. But we already saw that, by later testimonia, `static' and `dynamic' should come down to the same [Simplicius [{\small
DK} {\scriptsize{29B 2}}]. Given that Plato and Aristotle had both access to
Zeno's book, how is this to be understood?\footnote{This is discussed in detail by G.
Vlastos in a highly illuminating article, ``Plato's testimony concerning Zeno of Elea'';
reprinted in  G. Vlastos, {\em Studies in Greek Philosophy. Vol. I: The Presocratics},
Princeton University Press, Princeton, New Jersey, 1993, pp. 264-300.} Especially since Plato equally stresses the underlying
unity of Zeno's arguments [{\em Parmenides}, 127(e)]: {\em ($\ldots$) Is the intention of your arguments to
vindicate, against all that is [commonly] said, that plurality does not exist? And you
think that each of your arguments is a proof of just that, so that you believe you have
produced {\bf \em as many proofs that plurality does not exist as are the arguments you have composed?}}\footnote{Vlastos's translation, VLAS, p. 273. My bold.} My impression is that Aristotle, giving up on Plato's line of attack, contents himself with a less ambitious claim which has at least the advantage that it can be physically grounded.\footnote{If this be true, than it explains why Aristotle redirected his attention from metaphysics to physics, and not the other way around, as it is generally held in the literature. I thank again F. De Haas for pointing out to me this consequence of my own analysis.} I think we touch the heart of the matter here. Indeed, it is clear even on superficial inspection that underlying Zeno's paradoxes some idea must be present common to them all.\footnote{This insight has been defended, but not conclusively shown by G.E.L. Owen, ``Zeno and the 
mathematicians'', {\em Proceedings of the Aristotelian Society}, {\bf 8}, 1957. In my paper on Zeno's paradoxes (see ft. 16), I demonstrated this unity rigourously.} Then whence do these different approaches come from?

\subsection*{\sc Plato's Conceptual Division }

Now let us look at how precisely this kind of ontological paradox informed the core of Plato's metaphysical set-up. In the {\em Phaedo}, the first text ever to contain a clear formulation of the principle of contradiction, there is an extremely illuminating discussion of the relationship between the joining of two physical parts and the
addition of two mathematical line segments, where Plato again refers to the ``contradictory'' viewpoints on the causes of plurality and change proposed by pre-Socratic philosophy.\footnote{Tannery refers in this context to the critique of Protagoras on the geometers as quoted in Aristotle [{\em Met.} B, ii, 997b(33)-998a(5)]: \emph{(...) for as sensible lines are not like those of which the geometrician speaks (since there is nothing sensible which is straight or curved in that sense; the circle touches the ruler not at a [single] point, as Protagoras used to say refuting the geometricians).} In a rare display of negligence, DK destroy
exactly the point of the argument by translating \grk kan<onoc \eng 
by \emph{Tangente} instead of \emph{ruler}!  See [{\small DK}
{\scriptsize{80B 7}}], transl. in vol. II, p. 266; P. Tannery, {\em o.c.}, pp. 396-397.} 
 Socrates complains that, since it allows for akin things to have contradictory
causes, while phaenomena clearly distinct become causally undistinguishable, its results
cannot be valid [{\em Phaedo}, 100(e)-101(a,b)]. He gives the example of the
difference between ``being two things'' and ``being a pair of things''; the latter a
formal, the former a physical fact. He also focuses on the relation between parts and
wholes in number theory, by comparing the generation of `two' out of `two ones' by
bringing them together, and the separation of `two' into `two ones' by dividing it [{\em
Phaedo}, 97(a,e)]. The Zenonian influence is plain, although scarcely
discussed in the literature. But there is an even more striking parallel. Plato's philosophy is the first to introduce the idea that {\bf \em every concept is itself some sort of extension, and thus can be divided.} This conceptual division moreover consitutes the core of the parallel between his theory of being and his theory of knowledge, and became known through Aristotle as {\em
diairesis}. It is, 	according to Plato, the technique shoring up properly practiced dialectics, and is extensively applied in the {\em Sophist}, the {\em
Statesman}, and the {\em Philebus}. We will return in more detail to the {\em Philebus},
where the link between conceptual and mathematical {\em diairesis} is established. As is well known. Plato in the {\em Sophist} defines `dialectics' as the art
of making the proper distinctions between the forms that instantiate themselves in and
through particular things [{\em Sophist}, 253d(1-3)]. In that dialogue, moreover, the
difference established between contraries and contradictions --- between praedicative and
existential paradoxes --- lays out the underlying ontological preconditions to which correct application of dialectics must comply.\footnote{{\em
Platonists who doubt that they are spectators of Being must settle for the knowledge that they are investigators of the verb `to be'.} G.E.L. Owen, ``Plato on Non-being'', {\em
Plato: a Collection of Critical Essays}, vol. i, G. Vlastos ed., Anchor/Doubleday, N.Y.,
1971, p. 223.} Now consider a concept as an extension, with a domain of applicability that starts `here' and ends `there'. How then are `here' and `there' to be determined? By the dissection of a concept into its constituting contraries says Plato, harking back to something already explained in the {\em Parmenides}, viz., the dichotomic way of reasoning. The method by which {\em diairesis} should be applied is then
demonstrated in the {\em Philebus} and the {\em Statesman} by means of examples. In the
latter dialogue, while trying to define the good statesman, Socrates and his friends find
out that the most long and cumbersome, but nevertheless the best way [{\em Statesman}, 265(a)] to discover the specific forms instantiated in a thing is by systematically
dividing the concept associated to it\footnote{The complex relation that exists within Plato's system between a thing and the name that it bears deserves a discussion of its own and will therefore remain out of our scope in the present paper. Suffice it to remark that the relation is neither purely conventional, nor a straightforwardly ontological attribution, leaving space for `false opinion' on the epistemological level, while the presence of sufficiently strong ontological ground for  `true knowledge'  by means of correct inference is established by means of the participation theory.} in opposing halves, like `living/non-living' [{\em id.}, 261(b)],
`feathered/unfeathered' [{\em id.}, 266(e)] or `odd/even', instead of arbitrarily
separating off a part --- `Greeks' vs. `barbarians'; `ten' vs. `all other numbers' [{\em
id.}, 262(d-e)], say. This process ends when one bumps on undetermined parts or {\em elements} or \selectlanguage{polutonikogreek} stoiqe~ia \selectlanguage{english}  [{\em Sophist}, 252(b3)], that are not themselves capable anymore of specifying underlying parts [{\em Statesman}, 263(b)]. The word specification makes sense here, because the element found can very well be a quantitative plurality. The number of steps needed to reach from the undetermined unity to this level of specification --- the proportion between part and whole --- then defines somehow the original
concept [{\em Philebus}, 16(d)].\footnote{I therefore only partically agree to Pelletier's analysis of what it means to ``carve nature at its joints''. It has a very precise methodological import, which allows one to distinguish between ``predicates that do not denote forms'', apart from ``dynamical'' considerations. F.J. Pelletier, {\em o.c.}, pp. xvi-xvii.} This, however, is often not possible, especially not
when the praedicates are relative properties like `warm/cold'; `short/tall' \&c. Their
opposites will run apart into infinity unless a limit, a  \grk p'erac \eng or boundary is put to
them, in order to find the good {\em measure} that guarantees their non-destructive
aequilibrium and thus the reality of the thing they describe [{\em Statesman},
283(d,e)]. Given the constraints inhaerent to our finite reasoning capacities, adequate criteria
need to be developed to put a {\em peras} to this infinity: {\em he
who begins with any individual unity should proceed from that, not to infinity, but to
a definite number, and (...) conversely, he who has to begin with infinity should not
jump to unity, but he should look about for some number representing a certain quantity
(...)} [{\em Philebus}, 18(a-b)]. Now in case the concept under consideration indicates
an essence, natural ways of division and limiting conditions are to a certain extend
available. When it is purely relative, like with e.g. the warm/cold-divide, it is much
more difficult to get hold of it and prevent division from running into infinity. In that case, a limit will have to be {\em imposed}, not completely arbitrary, but nevertheless from the
outside. It is a matter of finding the correct {\bf mean} [{\em Philebus}, 13(c,d)]. The
Stranger answers to the young Socrates's question on the correct ways of division that
there are two possibilities: {\em In this way: the one through their [the beings] shared
relative greatness and smallness, the other through the necessary essence of those that
become.\footnote{This part of the passage is my own translation. The translations
devoted to it in the standard literature excell by their incomprehensibility. The idea
behind it is that if merely relative, then greatness and smallness are divided by
imposing a more or less arbitrary measure, while when concerned with things that become, i.e., real entities, division has to follow their intrinsical nature. The translation
offered by De Win in his Dutch translation was helpful in my rendering of this part of
the passage, since it does capture precisely that idea. X. De Win, {\em Plato. Verzameld Werk},
Boek III, De Nederlandsche Boekhandel, Antwerpen, 1980 (3th ed.), p. 564. I used older edition, because the more recent edition has been mutilated zealots of linguistic ``modernisation".} (...) we must
suppose that the great and small exist and are discerned in both these ways, and not, as
we were saying before, only relatively to one another, but there must also be another
comparison of them with the mean} \, [{\em Statesman}, 283(d-e)]. In order to acquire
knowledge of a something, one has to disscociate this one into a many {\em via the great
and small} by chopping it at every level into praedicative halves (more precisely: logical
complements) which represent the same property equipped with a positive or a negative valuation indicating its being or being not so-and-so. The {\em locus classicus} for the method in the dialogues is a passage in the {\em Phaedrus}, where he compares the division of a concept in opposing
parts with the dissection of a physical body: {\em The second principle is that of
division into species according to the natural formation, where the joint is, not
breaking any part as a bad carver might. Just as our two discourses, alike assumed,
first of all, a single form of unreason; and then, as the body {\bf \em which from
being one becomes double and may be divided into a left side and right side, each having parts right and left of the same name} --- after this manner the speaker proceeded to divide the parts of the left side and did not desist until he found in them an evil or
left-handed love which he justly reviled; and the other discourse leading us to the
madness which lay on the right side, found another love, also having the same name (...)} [{\em Phaedrus}, 265(e)].  This ``process of division by genera'' [\grk t~wn gen~wn dia'iresic\eng] is criticised by Aristotle [{\em Anal. prior.} I, xxxi] for being impossible with regard to completeness,
and for being question-begging: it is not possible {\em to effect a demonstration of
substance [\grk o>us'ia\eng] and essence [\grk t`o t'i >estin\eng]} this way, he says, because since
many different divisions are possible, the chosen one will be either arbitrary, either
rest on foreknowledge precisely on the points the division was intended to expose. This is the zest of the so-called {\em Socratic fallacy} of which Geach and his followers accused Plato. But Plato explicitly says that he choses {\em dichotomic} division because it {\em reduces} arbitrariness, and gives the strongest guarantee
that all specific forms instantiated in the being investigated may be correctly identified [{\em Statesman},
262(b); 265(a)]. The process, moreover, can in principle always be carried through up to
infinity, a consequence of the One-Many nature of the relation between things and ideas: {\em (...) all the
things which are ever said to exist are sprung from one and many and have inherent in
them the finite [peras] and the infinite [apeiron]. This being the way in which things are
arranged, we must always assume that {\bf \em there is in every case one Idea of everything} and
must look for it --- for we shall find it [the Idea] is there --- and if we get a grasp
of this, we must look next for two, if there be two, and if not, for three or some other
number; and again we must treat each of those units in the same way, until we can see
that {\bf \em the original unit is one and many and the infinite}, but just how many it is. And we
must not apply the idea of infinite to plurality {\bf \em until we have a view of its
[plurality's] total number between infinity and one;} then, and not before, we may let
each unit of everything pass on unhindered into infinity.}\, [{\em Philebus}, 16(a-d)]. The point by now should be clear.  Of course infinity belongs itself to the
realm of the Ideas, while our imperfect reasonings are limited to the finite realm.
Therefore a {\em peras} or boundary has to
be imposed, in order for them for reach sensible conclusions instead of no conclusions at
all. And we saw that, in case the concept to be divided is not relative, but contains
an essence --- a forerunner to Aristotle's kategorein hos ousias --- division will not be purely arbitrary,
and carry on until an intrinsic {\em peras} is reached. Indeed one can say that the Form
embodied by something here-and-now literally delineates it from the rest of the world: it grants its identity, it acts as individuator. The problem the Sophists have is that they fail to see
this latter point: {\em But the wise men of our
time are either too quick or too slow, in conceiving plurality in unity. {\bf \em Having no
method, they make their one and many anyhow, and from unity pass at once to infinity};
the intermediate steps never occur to them. And this, I repeat, is what makes the
difference between the mere art of disputation and true dialectic} [{\em Philebus}, 16(a-d)]. Plato moreover is very precise on the limiting criteria appropriate to this
aim. In order for division to be not arbitrary, but in accordance with the specificity of the subject of investigation, it should be ``through the middle'' (\grk di`a m'eson\eng), i.e., following {\em the longest way}
[{\em Statesman}, 265(b)]. To divide the numbers by cutting of two thousand as
opposed to all the rest, or to divide humanity by cutting of the Greeks as opposed to
all the foreigners --- the Lydians, the Phrygians could do the same ---, is not only
awkward, but arrests at once the process of rational division. How then to divide further?
Dividing numbers into even and odd (not-even), and humans into male and female
(not-male) is valid for everyone, and leaves open a number of further possibilities
[{\em Statesman}, 262(c)-263(e)]. Again, setting aside simply a number of a group encompassed by a basic essential idea only makes sense after {\em diairesis} up to that basic idea has been completed. Division through the middle is the best guarantee to find back the relevant subclasses at every step: {\em
It is safer to proceed by cutting through the middle, and in that way one is more likely
to find classes} [{\em Statesman}, 262(b)]. Now in order to find the {\em smallest}
relevant part, the essential idea --- the species, {\em eidos} in the strict sense of the word --- encoded in the original concept, the division has to be executed up to the point where a part is reached that is not itself anymore a class --- not a \grk g'enoc \eng [{\em genos}].\footnote{With Plato this is not yet a precise technical term; we translate ``class'' following Balme:  {\em \grk g'enoc \eng used in general as a ``kind'' composed of related members, would be a natural choice for a class-word; it is indeed ready to mean the genus which is divisible into related species.} D.M. Balme, ``\grk g'enoc \eng and \grk e~>idoc \eng in Aristotle's Biology'', {\em Classical Quartely}, {\bf 12}, 1, 1962, pp. 81-98.} These two conditions are moreover connected: {\em We must not take a
single small part, and set it of against many large ones, nor disregard species in
making our division} [{\em Statesman}, 262(a-b)]. The relationship is not reciprocal: although a class always is a part [{\em meros}], a part need not be a class at all: {\em That when there is a class of anything, it must necessarily be a part of the thing of which it is said to be a class; but there is no necessity that a part be also a class} [{\em Statesman}, 263(b)]. Indeed, once the level of the species reached, only quantitative divisions remain possible, eventually up to infinity. This is why the Peripathetic School will define a {\em species} by an elliptic
formula as {\em quae de pluribus differentibus numero in eo quod quid
praedicatur.}\footnote{Porphyry, as cited in J. Gredt, {\em Elementa Philosophiae
Aristotelico-Thomisticae}, Sumptibus Herder, Barcinone, 1961, p. 136.} Real {\em
diairesis} has finished: the difference between the Greeks and the Lydians does not
permit essential discrimination because their difference does not cut through the whole
of humanity anymore --- Greek women and Lydian women both remain women\footnote{These
examples are discussed in I.M. Crombie, {\em An Examination of Plato's Doctrines. Volume
II: Plato on Knowledge and Reality}, Routledge, London, 1963, pp. 371-372.} ---; it is,
whatever its importance, contingent, or as Aristotle [{\em Met.} Z 4 1030a(10-13)] would say,
accidential. This fundamental idea, that division discriminates between species up to
the point where the smallest difference that is not merely numerical is captured, lies
at the origin of the Peripathetic {\bf differentia specifica}: {\em Species infima (quae
sola est species in sensu stricto) dicitur etiam species specialissima, et differentia
infima, differentia specifica.}\footnote{J. Gredt, {\em Elementa}, p. 142.}\\

\par {\sc An Excursus on basic Number Theory}.--- In what follows some notions from number theory will be used; the next paragraph only serves to give the unfamiliar reader an intuitive grasp of the concepts involved. Consider a set of objects sharing some common property. How can we characterise this set without referring to its specific objects or properties? Clearly by some quantitative procedure: I can decide that one set is bigger than another one by counting their respective elements. But then the question arises what `counting' really means: it means to pair univocally every element of a given set to the elements of the set of natural numbers $\nn$, in their standard order. This is easy enough when the number of elements is limited. However, when we look at it carefully, we see that two different characteristic numbers appear: {\bf ordinal} and {\bf cardinal numbers}. Cardinals simply refer to the quantitative range of a set seen as a totality, while ordinals take the {\em order of appearance} of the elements into account. If I count a basket filled with apples, I can state that there are ten of them after completing the process, or I can say that I have the tenth apple in my hand, thereby stating that the ordinality of the basket is ten as well, although there also is an eighth apple and a ninth. In the case of finite numbers of elements, this all remains fairly straightforward. The trouble begins when infinities are involved. Indeed, although it is impossible to count an infinity, it is not impossible to conceive of one as a whole. This is what we do when we speak of ``the real line'', or any other geometrical object. Let us look at a simple example. Suppose we count the set of even numbers $\{ 2, 4, 6, 8, \ldots\}$, by associating to them the naturals in their standard order: $f(1) = 2, \,f(2) = 4, \,f(3) = 6, \ldots$ Although we would not expect it, the {\em countable} infinite cardinality (called Aleph-null, $\aleph_0$) of the two sets is the same, for I will always find a natural number to mark out whatever even number. So the part is as big as the whole, so to speak. Their infinite ordinality $\omega$ is also the same, because I count them in the same order. Now count the set of even and odd numbers, organised in sequence one after the other: $\{ 2, 4, 6, 8, \ldots, 1, 3, 5, 7, \ldots\}$. Their cardinality is again $\aleph_0$, but the latter's ordinality is 2$\omega$! Even worse, cardinalities greater than $\aleph_0$ do also exist. Let us look once more at a finite example. The set $\{1, 2, 3\}$ can be divided in all its parts or {\em subsets} by systematically labelling all elements for all possible subsets with a membershipsrelation `yes' or `no': $\{\varnothing, \{1\}, \{2\}, \{3\}, \{1,2\}, \{1,3\}, \{2,3\}, \{1,2,3\}\}$. We call the result the {\em powerset} of the original set. How many elements are there in this set? The number is always bigger than the number of elements $n$ of the original set, and  is given by $2^n$; in this case $2^3 = 8$. The cardinality of the powerset of $\nn$ is $2^{\aleph_0} = \aleph_1$; it is an {\em uncountable} infinity, {\em bigger} than $\aleph_0$. So infinity comes in kinds. In fact, this process can be continued, so that there are arbitrary many infinite cardinalities out there... This shocking fact was discovered by G. Cantor towards the end of the nineteenth century.\footnote{G. Cantor, ``Grundlagen einer allgemeinen Mannigfaltigkeitslehre'', in: {\em Gesammelte  Abhandlungen}, Georg Olms Verlag, Hildesheim, 1962. For a more reader-friendly treatment, see R. Rucker, {\em Infinity and the Mind. The Science and
Philosophy of the Infinite}, Princeton University Press, Princeton, 1995 [1982].} And as we have seen: one of the ways to get there is by a modern variant of dichotomic division. This will prove relevant for the remainder of our story.\\

In the paragraph of {\em Metaphysics} A [987b25-27] where Aristotle explains how knowledge and Forms hang together, he introduces seemingly out of the blue the concept of number, and comments that Plato disagreed with the Pythagoreans in making the in(de)finite {\em apeiron} not One but something he labels Indefinite Twoness. Put otherwise, he says, the infinite is made up from the Large-and-Small. Let us investigate in more depth the nature of the `twoness' incorporated in Platonic dichotomy. It is, as Aristotle explains in Book
$N$ of the {\em Metaphysics}, closely related to Platonic number
theory. The topic has caused considerable contention in the literature, stirring a debate on plato's so-called Unwritten Doctrine. I do not need to get involved, however, for I subscribe fully and completely to Sayre's point of view that {\em the main tenets attributed to Plato in the first book of the Metaphysics in fact are present in the Philebus}.\footnote{This fact is in itself quite undeniable and its recognition does not commit one to a position on the ``unwritten doctrine'' central to the quarrel surrounding the ``T\"{u}binger Schule''-interpretation (Kr\"{a}mer, Gaiser)  of Plato's philosophy. It should be remarked, though, that traces of Platonic {\em diairesis} can be found back in the earliest dialogues, and arguably contributed to the development of Plato's theory of Forms, at least according to M.K. Krizan, ``A Defense of {\em Diairesis} in Plato's {\em Gorgias}, 463e5-466a3'', {\em Philosophical Inquiry}, XII(1-2), 1-21. For the older scholarship on the question, see L. Robin, {\em La th\'{e}orie platonicienne des
id\'{e}es et des nombres d'apr\`{e}s Aristote.
\'{E}tude Historique et Critique}, Georg Olms Verlag [reprint 1908], Hildesheim, 1963; and J. Stenzel (see further). Contra Robin, Stenzel,  see H. Cherniss, {\em The Riddle of the Early Academy}, University of California Press, Berkeley, 1945. For a recent and moderate overview of the issues at stake,
see D. Pesce, {\em Il Platone di Tubinga, E due studi sullo Stoicismo}, Paideia, Brescia, 1990.  } We will follow Stenzel's in depth analysis\footnote{I agree to a certain extend with M. Dixsaut, {\em o.c.}, p. 221. sq. criticisms on Stenzel's version of Platonic diairesis in his
{\em Studien zur Entwicklung der Platonischen Dialectik von Sokrates zu Aristoteles}, Teubner Leibzig und Berlin, 1931. I do believe, however, that these objections disappear  when one takes Stenzel's {\em Zahl und Gestalt Bei Platon und Aristoteles}, Teubner, Leipzig, 1933 into account.}, and add to it a further elaboration and some illuminating examples provided by H. Oosthout in his commented translation of
{\em Metaphysics} $N$. We discussed before that, according to Plato, someone searching for knowledge about a given thing has to lay bare its
consituting ele\-ments by determining by means of dichotomic division the number of pairwisely opposite intermediate stages that separate (``are between'') that concept as a `one' and as an indeterminate [{\em apeiron}] `many' of characteristics --- \grk metax`u to~u >ape'irou te ka`i to~u <en'oc \eng [{\em Philebus}, 16(d)]\footnote{The procedure is marked out by
Plato's method of separation of the different kinds of `to be' as outlined in the {\em Sophist}: recognition of the intermediaries as praedicative contraries at once prevents one from falling into (existential)
contradiction.} ---, The intellectual activities concommittant are the work of the true
philosopher vs. the work of the sophist [{\em Sophist}, 253(c)].\footnote{J. Stenzel,
{\em Zahl und Gestalt Bei Platon und Aristotles}, Teubner, Leipzig, 1933 (2nd ed.), pp.
12-13.} It is obvious that Aristotle's later definition of dialectics finds its origin in
this distinction. In the divisional scheme proposed by Stenzel-Oosthout\footnote{H. Oosthout, {\em
Aristoteles. De Getallen en de dingen. De Boeken $M$ en $N$ van de Metafysica. Deel 2:
Boek $N$}, Klement, Kampen, 2004, pp. 74-79.} the opposing motions\footnote{Te dynamics
inhaerent in this time-dependent process has been analysed brilliantly by G. Deleuze,
while considering the paradoxical growth of Lewiss Caroll's {\em Alice in Wonderland}. See
G. Deleuze, {\em Logique du sens}, Minuit, Paris, 1969, pp. 9-12.} inhaerent in the
Large-and-Small come about through {\em aoristos duas} with a conceptual unity --- 1 --- at its centre: doubling in each step gives rise to the powers of 2 over the totality of the process, while halving in each step gives us the powers of 2$^{-1}$. Two modes of
generation are therefore at work simultaneously in the Platonic tree: doubling (2$^n$) and
halving (2$^{-n}$) of the original One:, $\ldots 1/8, 1/4, 1/2, 1, 2, 4, 8,
\ldots$.\footnote{H. Oosthout, {\em o.c.}, p. 79. Thus it can be linked to modern developments in the foundations of mathematics.  Hermann Weyl stresses exactly this aspect of Brouwer's intuitionistic approach, and links it explicitly to Plato: {\em Brouwer erblickt genau wie Plato in der Zwei-Einigkeit die Wurzel des mathematischen Denkens.}H. Weyl, {\em Philosophie der Mathematik und Naturwissenschaft}, M\"{u}nchen, 1926, pp. 43-51.} The process stops, says Plato, when one reaches the stage where one finds a part that is itself not further
divisible, i.e., which is itself not a {\em genos} of something else [{\em Statesman},
263(a,b)].	In modern terms: when one hits on those members or subsets  --- he calls them
{\em stoicheia}, elements --- which contain themselves not subsets of the originally given
set anymore, i.e., by filtering out the singletons from the powerset. In order to be exhaustive, this
means that you have to take all the possible subsets of a given starting set --- {\em genos} ---
by a yes/no-divisional procedure, exactly in the way a present-day mathematician would
proceed when constructing the powerset. This also explains why Plato insists on dichotomic division ``through the middle'' as
the method of ``the longest way'': it is the only one that guarantees that you will find
back the  truly ultimate members and not just some arbitrary intermediate result [{\em
id.}, 265(a-b]. {\em Aoristos duas} can therefore be understood at each stage as  incorporating $\{0, 1\}$, the
set of outcomes of the as yet undecided membershipsrelation. Evaluating all possibilities
comes down to taking
$\{0, 1\}^n$, with
$n$ the number of elements of the original set.\footnote{The relational logic that
underpins Plato's participation theory can be seen from this light. It also sheds more
light on the precise nature of the difference between Plato's and Aristotle's logical
systems.} The procedure is akin to what a botanist does when he or she consults a flora
and decides systematically whether a certain property is applicable to a given specimen
or not, thus narrowing down the number of remaining possibilities. I think it therefore appropriate  to retain a possibility which was considered but then dismissed by Sayre as a valid representation of Plato's procedure: {\em to conceive unity and being as states in a binary code (...), with succesive divisions of constituents ordered by level in an ``inverse tree'' (...).}\footnote{K.M. Sayre, {\em o.c.}, p; 55.} This idea gets confirmation from an unexpected side. In his notorious book  {\em Philosophie der Mathematik und Naturwissenschaft}, Hermann Weyl discusses L.E.J. Brouwer's intuitionistic approach to number theory, and links it explicitly to Plato: {\em Brouwer erblickt genau wie Plato in der Zwei-Einigkeit die Wurzel des mathematischen Denkens}.\footnote{H. Weyl, {\em Philosophie der Mathematik und Naturwissenschaft}, M\"{u}nchen, 1926, pp. 43-51.} Weyl refers to Plato's \grk >a'oristos duas \eng [indeterminate twoness], the principle by which {\em diairesis} [conceptual division]  generates both the large and the small. Every branch of the divisional tree represents a natural number which encodes for
a possible combination of properties, and, as Oosthout rightly observes, it would be the
fully Platonic option to consider the branch with a
number build up by only ones as the only one properly determined. I therefore propose
as a name for this Platonic quantification of {\em metaxu} well-definedness {\bf $\phi$-number} (philosophical number). Consider the example below. The number of
well-determination in the example in the scheme is at the bottom of the most right
branch of the divisional  tree, viz., 1111 or decimal 15. It is an ordinal number (reached stepwise), and contains full
information concerning the divisional pathway followed to define the element, for it
orders the determining sequence by means of {\em first differences}, i.e.
lexicographical.\footnote{K. Kuratowski
and M. Mostowski, {\em Set Theory}, N. Holland, Amsterdam, 1968, pp. 224-227.} The number of intermediate steps or generations needed to reach the bottom level in our example is 3, whence another relevant number arises as a consequence of the doubling-part of the
diairetic procedure, equal to
$2^3 = 8$, the number of {\em elements} found on the vertical line, or the cardinality of the subset
to which our $\phi$-number belongs. This is what Stenzel intended when he writes: {\em
Denkt man an die alte Darstellung der Zahlen durch punkte und fa{\ss}t man, was sehr nahe
liegt, innerhalb der Zahlengestalt jeden punkt als stelle auf, so ergeben sich {\em mit
einem Schlage Kardinal- und Ordinalzahlen,} und es zeigt sich sofort das einfache Bild
einer entstehung der Zahlenreihe durch stete anwendung der Zweiheit auf die eins und jede
sich ergebende neue Einheit in ihrer ''Zweifachmachenden'' (\grk diqotom'ac\eng) Natur.}\footnote{J. Stenzel, {\em o.c.}, pp. 30-32.}
As long as the number of steps needed
to reach the elements remains finite, the cardinality and ordinality of the sets of
products of the halving and of the doubling parts of {\em diairesis} will be the same.
This aequivalence between cardinality and ordinality is generally true in the finite
case. I accept with Stenzel-Oosthout that this is a modern formulation of what Plato had
in mind with his ---limited, i.e., finite --- Idea-Numbers. This can be represented graphically (see figure). The asymmetry
between 0 ({\em aoristos duas}) and 1 at every node of the diaretic tree is exemplified in the undeterminatedness of the
0-side at every bifurcation. And we saw that for Plato every chain of reasoning accessible to
the human intellect is by necessity finite. But...

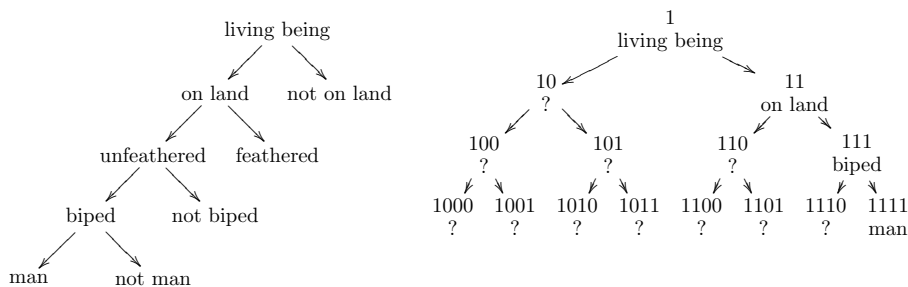
\begin{figure}[!ht]
\[
\resizebox{\textwidth}{!}{{\xymatrix@!0@R=3em@C=3em{&&&& \txt{living
being} \ar[ld] \ar[rd]\\
&&&\txt{on land} \ar[ld] \ar[rd] && \txt{not on land}\\
&&\txt{unfeathered} \ar[ld] \ar[rd] && \txt{feathered}\\
&\txt{biped} \ar[ld] \ar[rd] && \txt{not biped}\\
\txt{man} && \txt{not man}}}
\quad
{\xymatrix@!0@R=3em@C=1.5em{&&&&&&& \txt{1 \\ living being}
\ar[lllld] \ar[rrrrd]\\
&&&\txt{10 \\ ?} \ar[lld] \ar[rrd] &&&&&&&& \txt{11 \\ on land} \ar[lld]
\ar[rrd] \\
&\txt{100 \\ ?} \ar[ld] \ar[rd] &&&& \txt{101 \\ ?} \ar[ld] \ar[rd] &&&&
\txt{110 \\ ?} \ar[ld] \ar[rd] &&&& \txt{111 \\ biped} \ar[ld] \ar[rd]\\
\txt{1000 \\ ?} && \txt{1001 \\ ?} && \txt{1010 \\ ?} && \txt{1011 \\ ?}
&& \txt{1100 \\ ?} && \txt{1101 \\ ?} && \txt{1110 \\ ?} && \txt{1111 \\
man}}}}
\] 
\caption{A formal representation of Plato's divisional procedure}\label{}
\end{figure}

\noindent But something strange happens when this divisional process is carried through up to infinity. Every single being at every moment instantiates many properties, it ``necessarily and always is becoming two'' [{\em Parmenides}, 142(e)7-143(a)1], so that the divisional process in principle never stops, although by imposing order it somehow limits the \grk >'apeiron pl~hjoc \eng [infinite multitude] [{\em Parmenides}, 143(a)2].\footnote{See Sayre, {\em o.c.}, pp. 54-57 for a discussion of the relevant passages of the {\em Parmenides}.} Now at generation $\omega$, the smallest order imposed by ``counting'' the set of natural numbers when the smallest $\phi$-number written as an infinite sequence of ones would be reached\footnote{Of course an impossibility for mortal men, for one cannot stepwise reach $\omega$, although the ordinality of the set of natural numbers exists as much as the natural numbers themselves...}, two {\em different} kinds of infinite cardinalities arise for the halving and the doubling procedure, namely the countable and the uncountable infinite. Stenzel-Oosthout do not seem to realise the cardinal difference between these,
for they seem to think that even in the infinite case the number of elements of the
double-halve sequence remains the same at both sides of the 1--axis. The same holds true for Sayre, who says: {\em The result would be equivalent to an expression of the natural numer series in binary notation}.\footnote{K.M. Sayre, {\em o.c.}, p. 56.} But this is wrong! Aristotle not only gives clarification to what we already know from the dialogues {\em Sophist}, {\em Statesman} and {\em Philebus}, but plainly states that Plato ultimately had {\em two} infinities in mind: {\em He [Plato] supposes {\bf \em two infinities}, the Large and the Small} [{\em Phys.} 203a 15]. With each of these infinities goes a different mode of realisation: {\em everything is infinite, either through addition [i.e. stepwise], either through division [i.e. simultaneously] (...)} [{\em Phys.} 204a 6].\footnote{See the discussion in J. Stenzel, {\em o.c.}, p. 30 sq.; p. 60 sq.} These two represent stepwise (halving; think of the standard views on Zeno's paradox) vs. simultaneous (doubling) division, or in the traditional vocabulary, the {\em potential} and the {\em actual} infinite. This articulates with another mistake: contrary to what our authors believe, it is doubling which gives rise to the Small, and halving which gives rise to the Large once infinity is reached. This had been proven already by Zeno. He showed that it is indeed the stepwise halving process that generates infinitely many extended parts
which are included --- like atoms --- in the whole and therefore can still be further divided, while the simultaneous totality of the doubling-procedure generates {\em another} infinity of unextended parts,
members of the whole without being included in it --- like mathematical points.\footnote{Cfr. W.E. Abraham, {\em o.c.}, p. 48.} Both --- the {\em megala kai mikra} --- are present in every finitely extended thing.
The Platonic parallel with Zeno's analysis will by now be clear: Plato's Large-and-Small \grk t`o m'ega ka`i t`o mikr'on\eng\footnote{Cited by Aristotle in
e.g. [{\em Met. A}, 987b(20)]. Cfr. J. Stenzel, {\em o.c.}, p. 6.} and Zeno's
large[s]-and-small[s] \grk m'egala ka`i mikr`a \eng [{\small DK} {\scriptsize{29B 1}}] refer both to an infinite and through and through division, the first to its operative principle, the latter to its
paradoxical result.\footnote{Although I
found it nowhere discussed in the literature. The singular in
Plato's expression stems from the fact that he applies the paradoxical property to the abstract principle by which division is obtained, while with Zeno the plural directly refers to its result.}  Plato's metaphysically marked preoccupation with concsistent --- i.e., finite --- pluralities forces him to recover the paradoxical infinity resulting from Zeno's deictic point of view by means of an unlimited
number of partial, {\em finite} perspectives that cover the totality of the divisional tree. They are represented by the tree's branches and are uncountable infinite in number. These partial perspectives are what is accesible to our human `mind's eye', while the simultaneous realisation of all possibilities is confined to the world of Ideas, and is accesible only through the Form `be', present in and through anything that `is' in whatever sense. Thanks
to an idea of M. Serfati, it will be possible to formally distinguish the set of
perspectives ({\em l'ensemble des points de vue}) from the tree itself.\footnote{M.
Serfati, ``Quasi-ensembles d'ordre $r$ et approximations de r\'{e}partitions
ordonn\'{e}es'', {\em Math. Inf. Sci. hum.}, {\bf 143}, 1998, pp. 5-26. The underlying structure is a so-called Post algebra, of which Priest's paraconsisten logics are a special case. See the next paragraph.} Everything is embedded by an {\em apeiron pl\={e}thos} network of relations to everything else in the world, so that, in order to truly know something you should be able to assign it its place in the totality of the diairetic tree, which means that you would have to overview it completely at once. But, according to Plato, the Form `be' is not in the same ontological realm as the beings themselves, therefore this absolute knowledge is not within human reach. The best one can strive for is partial knowledge against the background of an understanding of how true knowledge in principle comes about, maybe after a walk outside the Cave. This sheds some light on the interconnection of some problems treated separately in the existent literature. At first the so-called ``Socratic fallacy'' imputed on Plato by Geach which we mentioned at the beginning of this paper. Geach's ``Assumption A'', that in order to {\em correctly predicate a given term `T', you must ``know what it is to be T'' in the sense of being able to give a general criterion for ta thing's being T} clearly is correct, for, metaphysically speaking, once you have to do it generally, you have to do it universally. And his Assumption B, that {\em it is no use to try to arrive at the menaing of `T' by giving examples} trivially follows from A, since it is impossible to pass in thought through infinitely many things --- as Aristotle takes care to remind us in [{\em Anal. Post.} I. 22, 83b7-8] --- especially when their number is uncountable. This by no means implies an error on Plato's side, given the careful discrimination of the ontological and eistemological levels involved. Secondly, this also clarifies an issue raised by Aristotle in [{\em Anal. Post.} II.13, 97a6-11], in a chapter that deals with {\em diairesis} and the  discovery of definitions, about ``some'' who claim that, on order to know the differentia of A, one has to know the differentia of B, C, D,... It is generally assumed that this argument undermines the possibility of giving any definition.\footnote{A. Falcon,  {\em Aristotle, Speusippus, and the Method of Division}, Classical Quartely, {\bf 50}, 2, 2000, pp. 402-414.} Tradition ascribes this view to Speusippus, Plato's student and successor in the Academy.\footnote{L. Tar\'{a}n, {\em A Critical Edition with a Collection of the Related Texts and Commentary}, Brill, Leiden, 1981.} I am not sufficiently familiar with the extant sources concerning Speusippus, but it seems clear to me that it is {\em not} recommendable to deny, as Falcon does, that Speusippus builds his system an absolute or universal diairetic tree. On the contrary, if he in any sense inherits the core structure of Plato's metaphysics --- and it seems to be the case that he does --- he needs it. Cherniss's viewpoint and mine are in agreement on this point, at least as I infer his from his description of speusippus endeavour: {\em For Speusippus, however, the essential nature of each thing is identical with the complex of all its relations to all other things, so that the content of existence is nothing but the whole network of relations itself, plotted out in a universal diaretical scheme. Consequently he could regard any particular being as the analogue of the point, for the different entities are simply different foci of the single system of relations. It is in the light of this theory of essence that the separately existing numbers are to be understood}.\footnote{H. Cherniss, {\em o.c.}, p. 42.} 

In the later tradition this will become the realm of the omnipresent God, who has all infinitely many possibilities of Creation at once lais out before Him.\footnote{I owe to a discussion with Michel Serfati that this basic scheme persists in Leibniz's philosophy. I am convinced that this is true as well for Newton, be it in a very different way.} For the human
mind, again, the relocation of the infinite to the realm of God's Ideas makes that the
part of the tree accessible to it always remains partial. This finite world is the world in which classical logic is valid. The paradoxical realm where all possibilites are realised at once is safely relocated to the Ideal World. With the pre-Socratics, however, the totality of the divisional process is paradoxically given in the tangible here-and-now. But then again with Plato, being ideal does not mean being unreal, quite the contrary! That is why the Platonic spectator of Being in its One/Many-marked appearance retains a definitely paradoxical flavour. The
differences between the Zenonian divider into {\em megala kai mikra} and the Platonic
spectator of {\em mega kai mikron} are nevertheless relevant, and can be clarified by
comparing the {\bf divisional ontologies} they represent.\footnote{When it comes to application of our inferences to reality, {\em A system of logic is a formalization of an ontology!} G. G\"{u}nther, ``Cybernetic Ontology and Transjunctional Operations'', BCL publication 68. Photomechanically reproduced from {\em Self-organizing Systems}, 1962, Yovits, Jacobi and Goldstein Eds., Washington D.C., Spartan Books, 1962, pp. 313-392.} By formalising them, these ontologies will become more transparant and more easily comparable. This is the subject of our next paragraph.

%\newpage

\subsection*{\sc Plato's not so classical classical logic}

In this section we shall study Platonic {\em diairesis} from the point of view of formal logic. The subject of Plato's logic is not a popular
one, although it deserves attention. The relational logic implicit in the participation
theory is, when discussed at all, often dismissed as being fallacious due to elementary
logical errors such as an improper `relation'-concept or self-praedication.\footnote{See e.g. C. Strang, ``Plato and the Third Man'', in: G. Vlastos, {\em Plato: a Collection of Critical Essays}, Doubleday, 1971, pp. 184-200.} However, it is
to be feared that the misconceptions involved are not Plato's. According to
H.-N. Casta\~{n}eda:
\begin{quote} {\em \small Contrary to the monolithic consensus among Plato scholars, in the
Phaedo Plato did distinguish, and soundly, between relations and qualities, and dealt
with genuine puzzles that arise in attempting to understand the nature of relational
facts. The reason why Plato's theory of relations has hitherto remained hidden to his
commentators is this: his commentators have either not understood the nature of
relations, or, more recently, they have adopted the dogma that a primary or simple
relation is just one atomic or indivisible entity that generates facts by being
instantiated at once by an ordered
$n$-tuple. (...) It might be suggested at this juncture that a nominalist must,
nevertheless, distinguish between a thing {\em a} being longer than another thing {\em b}, and the former being heavier than the second, and this distinction must lie in facts
themselves, in nature.}\footnote{H. N. Casta\~{n}neda, ``Plato's theory of relations'', in {\em Exact Philosophy}, Mario Bunge ed., Reidel Publishing Company, Dordrecht, 1973.}\end{quote} Our perspective does not focus on the relational
nature of Plato's logic, but on its inferential structure, c.q., on conceptual {\em
diairesis}, although we will be able to contribute something to the former aspect toward the end of this paper. The general idea expressed in the quote remains valid for us as well, however: one cannot properly appreciate the value of Plato's logic if one does not take its ontological foundation into account.

We discussed before why and how Platonic {\em diairesis} aimes at conceptually fleshing out the different Forms or Ideas that instantiate themselves in a specific being: {\em For it is from the
mutual intertwinement of the Forms that reasoning [logos]\footnote{We should never forget
that `to reason' and `to speak' at that moment in history are still understood as being the same thing!
R.B. Onians, {\em o.c.}, p. 13. Cfr. [{\em Phaedrus}, 266(b)].} comes forth} [{\em
Sophist}, 259(e)]. {\bf Plato's logic is the epistemological face of his participation
theory.}  Apart from what we discussed already, we are not going to dwell on its details or its evolution throughout
the dialogues, or on the disagreements concerning these points among
scholars.\footnote{I nevertheless refer once again to Dixsaut's and Sayre's excellent work. An overview of the topic and further references can be found in M.K. Krizan's {\em Defense}-article.} The main aspects of the divisional procedure are outlined clearly
enough by Plato himself in what I call the {\em diairetic dialogues}, --- {\em Sophist, Statesman} and {\em Philebus}, with the {\em Parmenides} as a requirement --- and suffice to expose its underlying logic. We shall show that
in Plato's system every being, when not yet divided, is a {\bf paradoxical
One-and-Many} --- the  \grk >'apeiron pl~hjoc \eng discussed in [{\em Parmenides}, 143(a)2] --- which becomes after complete {\em diairesis} a {\bf paraconsistent
Large-and-Small.} Indeed, in [{\em Philebus}, 18(e)] Plato insists not only on the infinity of individual things, but also on the infinity {\em within} every individual thing. It
really is a principle, {\em (...) its nature is quite marvellous, for that one should be many or many one, are wonderful propositions (...)} [{\em Philebus}, 14(c)]. But when we leave it at that, paradoxes involving infinity will pop up
that cause problems for our thought, especially when we use general concepts: {\em
Those (...) are the common and acknowledged paradoxes about the one and many, which I may say that everybody has by this time agreed to dismiss (...)} [{\em Philebus}, 14(d)]\\

\par {\sc An Excursus on Paraconsistency}.--- Let us, before we give an outline of the formalism in which Plato's logic can be cast, present an overview of some basic notions, whereby comprehensibility will be given precedence over rigour. In general, a {\bf logic} is a system of propositions $P,Q\ldots$ and operators $\wedge, \vee, \neg,\dots$ which allow to connect them into wellformed formulae according to some set of rules specified in advance. A {\em proof} \/ is a deductive scheme by which one formula can be transformed into another, given certain hypotheses and by careful application of the given set of rules (note that this presupposes some notion of identity). There will be encoded in one way or another the notion of {\em entailment}, which captures the general idea of `consequence'. This is linked to the presence of a formally workable implication-relation. We can consider our formal system from two different points of view: the {\em syntactical} viewpoint concerns the correct construction of formulae from the basic formal elements, while the {\em semantical} point of view investigates the validity of an inference by means of the algebra of the truth values involved. We will take the latter point of view in what follows. A {\bf classical logic} (in the traditional sense) has as basic operators `AND', `OR' and `NEGATION', and rests for the validity of its inferences on three basic principles, the principle of contradiction: $\sim \! P \wedge \!\sim P$ [PC], the principle of identity:  $P = P$ [PI], and the Excluded Middle or Tertium non datur:  $P \vee \sim P$ [TND].\footnote{It is common practice to formally present classical logic without making any reference to them, but it suffices to look at the truth table values for the material implication to realise that they encode the {\em ex falso}, hence the PC. I realised this during a discussion with Koen Lefever on a system he proposes to classify logical systems.} They can already be found --- in implicit or explicit form --- in Aristotle's account in the metaphysics and in the books of the {\em Organon}.\footnote{Aristotle did never formulate explicitly the PI. He did formulate the TND. However, it is the PC that is introduced as an axiom on ontological grounds, and translated in epistemological terms: {\em Metaphysics} \grk G\eng, 1005a(19)-1005b(33).  This is exactly why it is logically unprovable: {\em Met.} \grk G\eng, 1005b(3)-1006(11). The TND is provable in propositional logic when one accepts the PC; it is therefore clear that the original constitutive priority lies with the latter.} They guarantee that propositions be either true or false, and nothing else. Classical logics therefore are two-valued: the set of truth values for a proposition is the pair $\{0, 1\}$. Such logics are also very appropriately called {\em explosive}, because even the smallest violation of the PC causes a total collaps of the inferential structure, captured in the ominous formula {\em ex falso quod libet}. Now {\bf paraconsistent logics} are inferential systems that accept {\em to a certain extend} violations of the PC, without automatically becoming explosive: {\em (...) a minimal condition for a suitable inference relation in this context is that it not be explosive. Such inference relationships (and the logics that have them) have come to be called {\bf \em paraconsistent}}.\footnote{G. Priest, ``Paraconsistent Logic'', {\em
Handbook of Philosophical Logic}, 2nd ed., Kluwer, Dordrecht etc., 2002, p. 288.}  There can be many good reasons for developing less constrictive
logical systems, where violation of the PC has more moderate consequences. Such non-classical systems have been studied throughout the
twentieth century, following the seminal work of Jan \L ukaciewicz. He critisised the
Principle of Contradiction, and its corrolary, the Excluded Middle, for basically
ontological reasons, and introduced the notion of many-valued logics instead.\footnote{J. \L ukasiewicz, {\em \"{U}ber den Satz des Widerspruchs bei Aristoteles},
trans. J.M. Bochenski, Georg Olms Verlag, Hildesheim, 1993.} Four valued logical systems have been introduced by Da Costa and been further developed
mainly by Belnap\footnote{N. Belnap, ``A useful four-valued logic'', in: {\em Modern
Uses of Multiple-Valued Logic}, J.M. Dunn and G. Epstein (eds.), D. Reidel,
Dordrecht-Boston, 1977, pp. 8-37.} and Priest. The approach developed by Priest is the most relevant here: certain contradictions have to be allowed since under certain conditions paradoxical situations do occur. The slogan would then be that (some) paradoxes {\em really} exist. This view is sometimes called {\em dialethism}\footnote{G. Priest, {\em In Contradiction. A study of the transconsistent}, Martinus Nijhoff, 1987.} and it touches of course upon Zeno's paradoxes, which as we saw are at the heart of Plato's restrictive intervention with respect to the relations that are possible between reality and what we can say about it.

\bigskip Now to come back to our main topic: the  point of view we will take with respect to Platonic {\em diairesis} is the algebraic one. This is natural enough, since it translates our semantical preoccupations and it allows us to see Plato's diairetic tree as an inverted propositional semi-lattice. A
{\em proposition} will be a (finite or infinite) branch of the tree
starting from a bottom element, the given concept. Logically, each branch is the conjunction, lattice-theoretically the
meet $\bigwedge p_i$ of determined and indetermined praedicates along its descending nodes. Here we view a concept extensionally as a set of elements --- `beings' --- on which a given praedicate --- `Form' --- is applicable. The division of a concept will be represented formally as the division of a given concept set into a smaller property set and its exclusive complement, in an exhaustive way. The partial order $\order$ over the semi-lattice expresses in a natural way the idea of logical strength as the diairetic `sharing this property'-relationship.\footnote{S. Vickers, {\em Topology via Logic}, Cambridge University Press, Cambridge, 1989/1996, p. 13.} This obviously coincides with set-theoretic inclusion.  Platonic division, when completed, thus comes down formally to taking the powerset of the original conceptset. This, however, does {\em not} yet exhaust the structure that is present in the diairetic tree. The intensional
nature of Plato's logic emerges clearly in the path-dependance of the truth
values at every node.\footnote{Intensionality on behalf of Plato has been argued for
by J.M. Moravcsik, ``The Anatomy of Plato's Divisions'', in: E.N. Lee, A.P.D. Mourelatos
and R.M. Rorty (eds.), {\em Exegesis and Argument. Phronesis Supplementary Volume I},
N.Y., Humanities Press, 1973, pp. 324-348.} Indeed an additional order $\prec$ has to be imposed to catch the fact that the path towards a specific stage in the division --- 
i.e., a unique $\phi$-number --- is retained: at every level the $\phi$-numbers will
again be ordered, according to first differences. Such an order is called {\em
lexicographical}. This allows to clarify Serfati's suggestion that
one can regard a higher order division as ``meilleur''\footnote{M. Serfati, {\em o.c.},
p. 6}, and formally captures the classic idea that species can be distinguished by their
{\em differentia specifica}: the smallest difference which does not merely discriminate
numerically between individuals of the same kind. Thus the overall order is not partial,
but total: every branch will be uniquely determined by the totality of the procedure
leading to it. Mathematically speaking, this procedure gives rise to a semi-lattice with a total order $\prec$. This is illustrated in the diagram below:

\begin{figure}[!ht]
\[
\vcenter{\xymatrix@C=.5pc@-1pc{&& &&&&&&&&&&&&&&
*++[o][F-]{\bot}
\ar[dllllllll]_-{0} \ar@{.>}[drrrrrrrr]^{\mathbf{1}} \\ 
&&
&&&&&&  *++[o][F-]{0}
\ar[dllll]_-{0} \ar[drrrr]^{1} &&&&&&&&&&&&&&&& *++[o][F-]{1}
\ar@{.>}[dllll]_-{\mathbf{0}}
\ar[drrrr]^{1} \\ 
&& && *++[o][F-]{00}
\ar[dll]_-{0}
\ar[drr]^{1} &&&&&&&& *++[o][F-]{01} \ar[dll]_-{0} \ar[drr]^{1} &&&&&&&& *++[o][F-]{10}
\ar[dll]_-{0}
\ar@{.>}[drr]^{\mathbf{1}} &&&&&&&& *++[o][F-]{11} \ar[dll]_-{0}
\ar[drr]^{1}\\  
&& *++[o][F-]{000} \ar[dl]_-{0}
\ar[dr]^{1} &&&& *++[o][F-]{001}
\ar[dl]_-{0} \ar[dr]^{1} &&&& *++[o][F-]{010} \ar[dl]_-{0} \ar[dr]^{1} &&&&
*++[o][F-]{011} \ar[dl]_-{0}
\ar[dr]^{1} &&&& *++[o][F-]{100} \ar[dl]_-{0} \ar[dr]^{1} &&&& *++[o][F-]{101}
\ar@{.>}[dl]_-{\mathbf{0}} \ar[dr]^{1} &&&& *++[o][F-]{110} \ar[dl]_-{0}
\ar[dr]^{1} &&&& *++[o][F-]{111} \ar[dl]_-{0} \ar[dr]^{1}\\ 
&  {\vdots} && {\vdots} && {\vdots} && {\vdots} && {\vdots} && {\vdots}
&& {\vdots} &&  {\vdots} && {\vdots} && {\vdots} && {\vdots} && {\vdots} && {\vdots} &&
{\vdots} &&  {\vdots}&& {\vdots}}}\\
\]
\caption{The diairetic semi-lattice of finitary arguments}\label{}
\end{figure}
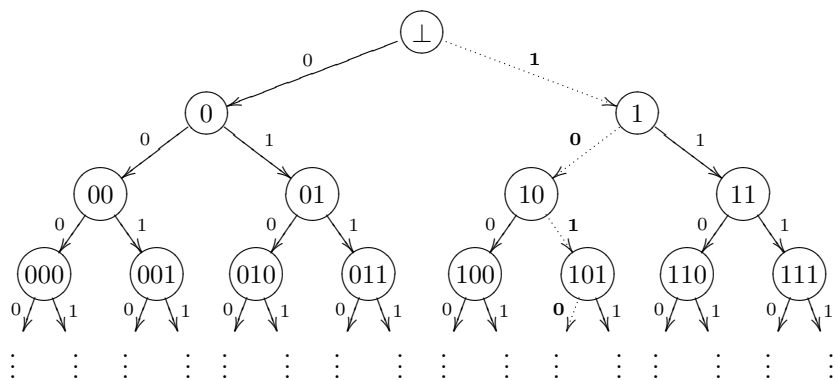

\noindent A {\em proof} of a
proposition consists of the whole diairetic tree up to the required generation. A
concept will be completely determined through {\em diairesis} up to generation
$p_n$ with $\nu(p) = 1, \forall \, n$. We can therefore speak of the determination value of every sub-branch taken as a whole. When only `$1$'s occur, this total value will be the $\phi$-number $<\! 1_n \!>$. Occurrence  of a `$0$' in a sequence $<\! p_n \!>$ means that there is still {\em aoristos duas} contained in the concept
under consideration, c.q., that determination up to that level is not logically
complete even though the appropriate singletons are reached, so that its value remains indetermined, a fact symbolised by zeros in the sequence: $<\!0/1_n\!>$. This is a way to diagnose that division has followed a wrong path. When division has been carried out correctly, i.e., through the middle, the level of the \grk >id'ia o>us'ia, \eng the esential ideas [{\em Phaedo}, 101(c)3)] instantiated in the being under investigation can be reached. If this is not possible, an appropriate measure or limit ({\em peras}) has to be imposed (we saw that perfect knowledge of all ideal relationships of a being would require knowledge of the uncountable Total Tree of Being, of which all finite trees are a part). The truth values in Plato's system will therefore be {\bf (partially/completely) determined} and {\bf (partially/completely) indetermined} rather than simply true and false.\footnote{I choose this terminology in reference to W.E. Johnson's discussion of ``the principles of logical division'' in his classic on logic, although my use of `determined' vs. `in[completely] determined' is different from his.W.E. Johnson, {\em Logic}, Part I, Chapter IX: ``the
determinable'', Cambridge, 1921.} Plato's {\em peras} criterion can thus be translated as the requirement that division, in order to be complete, should be carried
out up to the level of the singletons in the powerset: they are ``parts but not
classes'' in Plato's terms, species but not genera, which, although themselves sets,
cannot be written as the union of smaller subsets anymore, hence cannot be further divided. The choice of the concept to start with will influence the
length of the tree, a point that elicits Aristotle's negative comment. Filtering the
singletons out will give the fundamental elements or {\em stoicheia}.
The fact that this element-species still is a numerical plurality --- the species `man'
obviously comprises many individuals --- does not change anything to the matter, as we discussed before. But what is the Platonic determination value of a proposition obtained
after complete 'up to singletons) division? Clearly we need to include the diairetic tree associated to a specific division in its entirety in our formalism. Once division is completed, Plato obviously does not care anymore about the intermediate inferential steps. It is therefore natural to associate the determination value `$1$'
to the branch with the  $\phi$-number $<\! 1_n \!>$, and the value `$0$' to all the
others. The totality of
the divisional process is represented by the completed tree which contains both determined and indetermined elements, and therefore has itself the value `both'. Finally, evidently nothing has been determined
when no division is accomplished at all (`none'). Now we put completely determined = top element = $\{0, 1\}$; completely indetermined = bottom element = $\varnothing$. The number of Platonic determination values can thus again be limited to four:
{\bf nothing, determined, undetermined, everything.} Set-theoretically we have
$\{\varnothing, \{0\}, \{1\}, \{0, 1\}\}$. The set of overall Platonic truth values
is the powerset of the set of truth values for classical
logic (the pair $\{0, 1\}$)! it can be represented by a powerset lattice:

\begin{figure}[!ht]
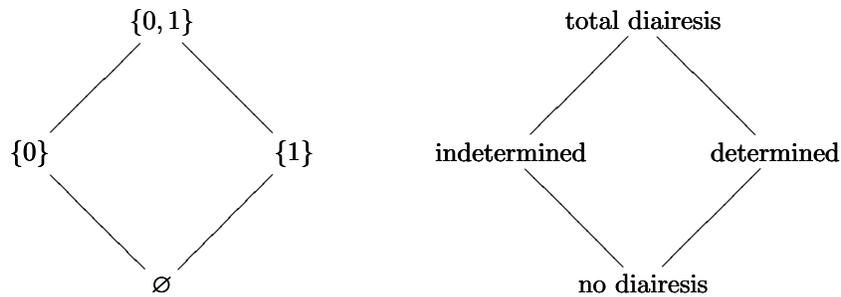

\[
\vcenter{\xy
\Atrianglepair(0,500)/-``-``/[{\{0, 1\}}`{\{0 \}}``{\{1 \}};````]
\Vtrianglepair/``-``-/[{\{0 \}}``{\{1 \}}`{\varnothing};````]
\endxy}
\qquad\qquad  
\vcenter{\xy
\Atrianglepair(0,500)/-``-``/[{\text{total
diairesis}}`{\text{indetermined}}``{\text{determined}};````]
\Vtrianglepair/``-``-/[{\text{indetermined}}``{\text{determined}}`{\text{no
diairesis}};````]
\endxy}
\]\caption{Plato's four valued truth-lattice}\label{}
\end{figure}

\newpage

\noindent The {\em determination value} of single propositions or branches is a mapping into the four valued set of possible determinations $\{\text{none, true, false, both}\}$. Every step in the argument adds a 0 (`false') or a 1 (`true') to the determination chain. The determination `both' means that, at the end of a reasoning process, one notes that one has a chain of either zeros, either ones, or both. The determination `none' means that no divisional resolution of the original {\em apeiron pl\={e}thos} has yet been executed: you see everything at once, but cannot discriminate anything at all. (In argumentative terms, to go from `none' directly to `both' would imply a forbidden un-limited One/Many-paradoxical jump.) In Platonic terms this comes down to a complete lack of knowledge exept with respect to the existence of the thing under consideration. We can make up the tables for the determination values of complex propostions and generalise them to all propositions in a given tree (finite or infinite):

\begin{figure}[!ht]
\[
\left.
\begin{array}{l|l}
			&	\neg	\\
\hline
t		&	f	\\
f		&	t	\\
b		&	b	\\
n		&	n		
\end{array}
\right.
\quad
\left.
\begin{array}{l|llll}
\wedge	&	t	&	f	&	b	&	n	\\
\hline
t		&	t	&	f	&	b	&	t	\\
f		&	f	&	f	&	f	&	f	\\
b		&	b	&	f	&	b	&	b	\\
n		&	t	&	f	&	b	&	n	
\end{array}
\right.
\quad
\left.
\begin{array}{l|llll}
\vee		&	t	&	f	&	b	&	n	\\
\hline
t		&	t	&	t	&	t	&	t	\\
f		&	t	&	f	&	b	&	f	\\
b		&	t	&	b	&	b	&	b	\\
n		&	t	&	f	&	b	&	n	
\end{array}
\right.
\]
\caption{Truth tables for Platonic {\em diairesis}}\label{}
\end{figure}

\noindent This is important, because the total number of possible $0/1$-combinations in the Total Tree {\em is} infinite. To be precise, the total number is $\{0,
1\}^{\nn}$. Recall that this does not require the number of inferential steps in a
specific argument to be infinite, on the contrary! Although the number of properties and
therefore the number of possible diairetic perspectives on every specific being is
infinite, a correct argument is kept  finite precisely by succesfully digging up the
appropriate {\em peras} or limiting condition, viz. the set of essential properties, to
the thing investigated.\\

\noindent Plato's methodological scheme thus has the structure of a {\bf four-valued paraconsistent logic}, although with the remarkable feature that it remains Boolean, because the complements are exclusive. An interesting point is furthermore that this logic can be reduced to first order
praedicate logic, as is the case for the relational system exposed by Casta\~{n}eda, so that our reconstruction of Plato's logical system remains at least in
principle commensurable to his results.\footnote{H. N. Casta\~{n}neda, ``Plato's theory of relations'', in {\em Exact
Philosophy}, Mario Bunge ed., Reidel Publishing Company, Dordrecht, 1973.} Plato's logic is not entirely classical, since it violates in a precise way his own dogma that sums up the demand for logical consistency, the Principle of Contradiction, by allowing somthing to be ``true'' and ``not true'' up to a certain extend, though without rendering the whole thing {\em explosive}. This lends some additional support to Sayre's point of view, that Plato shifts in the later dialogues to an ontology which is less restictive with respect to the self-identity of the Ideas. This relaxation never occurs on the level of single propostions; it implies consideration of the whole of which a being is part, or, epistemologically speaking, of the relational context of the proposition. Only on the level of the Total Tree this process would come to an end, since every being after finite division still contains some undetermined twoness, and therefore remains uncompletely resolved. On the level of the supreme idea of Being the undetermined twoness ultimately disappears. On the level of Being it would cause not innocent contraries, but fatal contradictions.\footnote{Prepared in the {\em Sophist}, where the different modi of the verb `to be' --- existential and praedicative --- are discussed in detail and related to their ontological import. Aristotle follows the lead in [Met. $\Delta$, vii], where he discerns the kinds of `be' not according to their function, but according to the kind of oppositions to which they give rise. For a discussion K. von Fritz, ``Der Ursprung der Aristotelischen Kategoriernlehre'', in: {\em Schriften zur griechischen Logik}, Fromann-Holzboog, Stuttgart-Bad Cannstatt, band 2, 1978, pp. 14-15.} Plato avoids this by his limiting or {\em peras}-condition which renders all particular arguments finite (we project them against a simplifying screen, so to say), while the actual infinity of possibilities in which they lie embedded captures the inconsistency required by the fact that real things are undetermined and therefore on a deeper level instable without {\em a priori} destroying the fabric of knowledge. 

\begin{figure}[!ht]
\resizebox{\textwidth}{!}{\includegraphics{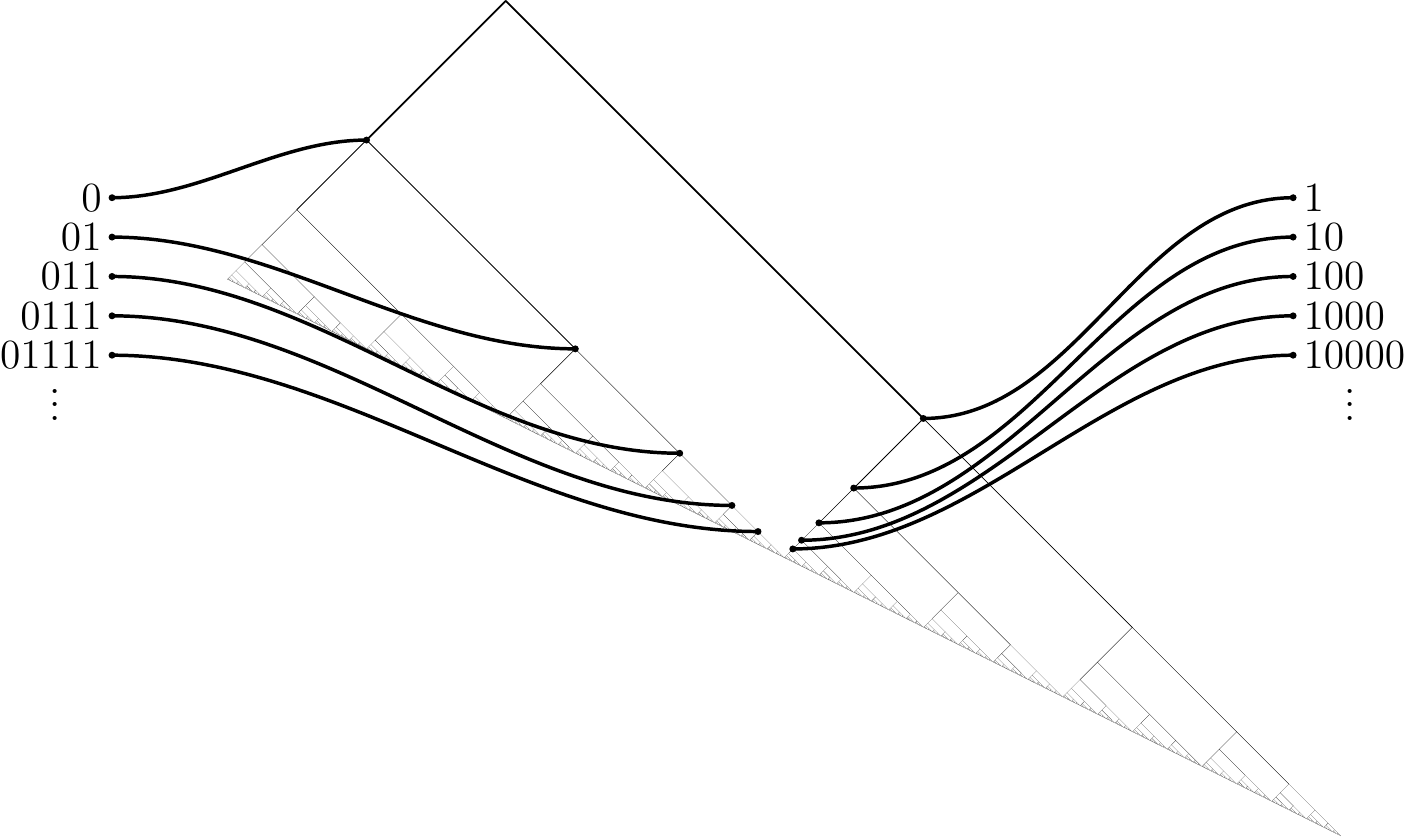}}
\vspace{-2cm}
\caption{The Philebian fractal}\label{Philebian_fractal}
\end{figure}

\noindent Plato's divisional logic can easely be generalised to an infinity of truth values, by taking the diairetic tree itself as the basis of a truth value semi-lattice with bottom element $\varnothing$ and ``top element'' $\{0, 1\}^{\small{\nn}}$.\footnote{M. Barnsley, {\em Fractals Everywhere}, Academic Press,
Boston Etc., 1988.} Every sequence encodes a numerical truth value between
$0$ and $1$; the generalisation thus naturally leads us into the realm of fuzzy logic.
One then construes the zeros in every determinative sequence as an explicit
quantification of conceptual uncertainty or indeterminedness, or as a statistical
distribution of the possible choices in a game. The geometry that goes with it is that of a  ``Theorem Fractal'' as proposed by Grim.\footnote{As was pointed out to me by J-P. Van Bendegem. See P. St. Denis and P. Grim, ``fractal Images
of Formal Systems'', {\em Journal of Philosophical Logic}, {\bf 25}, pp. 181-122, 1997.} When referring specificially to Plato's construction, I will call it the {\em Philebian fractal}. A modal
interpretation of the Philebian fractal, and thus of the Total diairetic Tree, could be that, taking the ideal world of full determination as the actual one,
the growing number of zeros in a formula at a given generation would indicate the
decrease in ``accessibility'' of the possible world encoded by it when compared to the one of the determined concept. This remains in line with our earlier discussion about {\em deixis}, because in Zenonian {\em this world} complete determination also coincides with total actualisation, but centered in the Here-Now. Modality for Zeno would be spatiotemporal disclocation from the pointing world-centre at the Here-Now. But in Plato's world, only for God everything is simultaneously present, as if only His eye can bear the sight of everything at once in all aspects of its being and non-being without being abhorred or disturbed by this unfathomable multiplicity. The differences between the Zenonian divider into {\em megala kai mikra} and the Platonic spectator of {\em mega kai mikron} are important, because transition form the One to the Many in the latter case can only take place in a mediated way. They remain connected nevertheless. That is why the Platonic spectator of Being in its One/Many-marked appearance definitely retains a paradoxical flavour. 

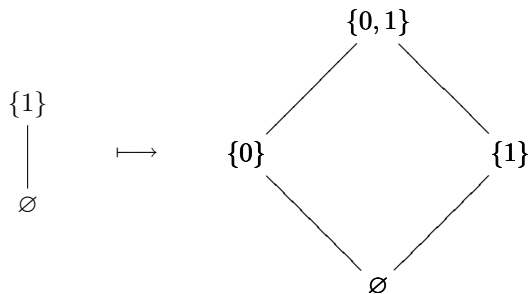
\begin{figure}[!ht]
\[
\vcenter{\xymatrix{{\{1 \}} \ar@{-}[d]\\ {\varnothing}}}
\qquad\longmapsto\qquad 
\vcenter{{\xy
\Atrianglepair(0,500)/-``-``/[{\{0, 1\}}`{\{0 \}}``{\{1 \}};````]
\Vtrianglepair/``-``-/[{\{0 \}}``{\{1 \}}`{\varnothing};````]
\endxy}}
\]
\caption{From Zeno's ``logic'' to Plato's logic}\label{Zeno's to Plato's logic}
\end{figure}

\noindent It are precisely these two infinities that will trouble Aristotle so deeply. For him it does not make sense to speak of a reality more real than reality as it is given. If another ontological level is needed, then it should be one that is {\em less} real, hence his relocation of ``actual'' infinity to the realm of imagination. It is not completely unreal, however, but has the status of a kind of limiting case, akin to the final cause of any concrete being. Aristotle realises that he needs actual infinity if he wants to have the continuum; that is why he does not simply can do away with it. Epistemologically, Plato's {\em peras} becomes Aristotle's {\em horos} in a syllogistic argument, classically known as the {\em terminus} or ``extreme'', of which there are two, the Major and the Minor, not incidentially the latin translations for ``Great'' and ``Small''. But that is a subject for another paper.

\bigskip
\subsection*{\sc Acknowledgements}
\smallskip

I wholeheartedly want to thank Frans de Haas (Leiden University) for many hours of fruitful discussions, and a wealth of useful suggestions and interesting ideas. During my stay at Leiden I had less time than I had wished to discuss the themes of my paper with Vasilis Politis (Trinity College, Dublin), who was visiting there in the same period. His influence on some of my work is nevertheless undeniable. I am also inspired by some discussions at the Leiden seminar on the {\em Analytica posteriora} 2007-2008 organised by Frans, and of which Vasilis and Pieter-Sjoerd Hasper (Groningen University) were part. Jean-Paul Van Bendegem (Brussels Free University) made some useful technical suggestions, and Tim Van der Linden (idem) the figures. Finally, I am much obliged to Claudio Beccari of the University of Turin, who not only gave me invaluable advise with respect to editing Greek texts in \LaTeX, but also took the pains to comment on the content of my paper. 

\bigskip
\subsection*{\sc A note on bibliography and used translations}
\smallskip
\par\noindent For Plato's text I relied on the Oxford edition by Burnett. The Aristotlian texts are from the Loeb edition. The pre-Socratic fragments all stem from the Diels-Kranz edition. The translations of Plato are based on the sources mentioned in the following list, but none of these are strictly followed and are adapted when deemed necessary: \\
\par\noindent
{\em Platonis Opera},  I. Burnet (ed.), Oxford Classical Texts, Clarendon, Oxford,1900--1907 [1952].
\par\noindent
{\em Plato, Collected Works in Twelve Volumes}, H.N. Fowler (ed.), Loeb Classical Library, Harvard University Press, Cambridge, Mass., 1914 [1999].
\par\noindent
{\em The Dialogues of Plato}, translated into English with Analyses and Introductions by B. Jowett,  in Five Volumes. Oxford University Press, 1892.   
\par\noindent
{\em Plato, Verzameld Werk}, X. De Win (vert.), De Nederlandsche Boekhandel/Ambo,
Antwerpen, 1980.                     
\par\noindent
{\em Aristotle, Complete Works in Twenty three Volumes}, Loeb Classical Library, Harvard University Press, Cambridge, Mass., 1930 [2002].
\par\noindent
H. Diels and W. Kranz, {\em Fragmente der Vorsokratiker}, erster Band,
Weidmann, Dublin, Z\"{u}rich, 1951 [1996].\\

\par\noindent For further bibliographical information the reader is referred to the
footnotes.

\end{document}